\documentclass{amsart}

\usepackage{slashbox}

\newtheorem{theorem}{Theorem}[section]
\newtheorem{lemma}[theorem]{Lemma}
\newtheorem{proposition}[theorem]{Proposition}

\theoremstyle{definition}

\theoremstyle{remark}
\newtheorem{remark}[theorem]{Remark}

\numberwithin{equation}{section}

\usepackage{color}

\begin{document}

\title[An integral involving the modified Struve function of the first kind]{Bounds for an integral involving the modified Struve function of the first kind}


\author{Robert E. Gaunt}
\address{Department of Mathematics, The University of Manchester, Oxford Road, Manchester M13 9PL, UK}
\curraddr{}
\email{robert.gaunt@manchester.ac.uk}
\thanks{The author is supported by a Dame Kathleen Ollerenshaw Research Fellowship.  }



\date{27 January 2021}

\dedicatory{}

\commby{}

\begin{abstract}Simple upper and lower bounds are established for the integral $\int_0^x\mathrm{e}^{-\beta t}t^\nu \mathbf{L}_\nu(t)\,\mathrm{d}t$, where $x>0$, $\nu>-1$, $0<\beta<1$ and $\mathbf{L}_\nu(x)$ is the modified Struve function of the first kind. These bounds complement and improve on existing results, through either sharper bounds or increased ranges of validity. In deriving our bounds, we obtain some monotonicity results and inequalities for products of the modified Struve function of the first kind and the modified Bessel function of the second kind $K_{\nu}(x)$, as well as a new bound for the ratio $\mathbf{L}_{\nu}(x)/\mathbf{L}_{\nu-1}(x)$.
\end{abstract}

\maketitle

\section{Introduction}\label{intro}

In a series of recent papers \cite{gaunt ineq1,gaunt ineq3,gaunt ineq8}, simple upper and lower bounds, involving the modified Bessel function of the first kind $I_\nu(x)$, were established for the integral
\begin{equation}\label{intbes}\int_0^x \mathrm{e}^{-\beta t} t^\nu I_\nu(t)\,\mathrm{d}t, 
\end{equation}
where $x>0$, $0\leq\beta<1$.  The conditions imposed on $\nu$ differed for several of the inequalities, although in all cases $\nu>-\frac{1}{2}$, which ensures that the integral exists.  For $0<\beta<1$ there does not exist a simple closed-form formula for this integral. The inequalities of \cite{gaunt ineq1,gaunt ineq3,gaunt ineq8} played a crucial role in the development of Stein's method \cite{chen,np12,stein} for variance-gamma approximation \cite{eichelsbacher, gaunt vg, gaunt vg2,gaunt vg3}. As the inequalities of \cite{gaunt ineq1,gaunt ineq3,gaunt ineq8} are simple and surprisingly accurate, they may also be useful in other problems involving modified Bessel functions; see, for example, \cite{bs09,baricz3} in which inequalities for modified Bessel functions of the first kind were used to derive tight bounds for the generalized Marcum $Q$-function, which arises in radar signal processing.

The  modified Struve function of the first kind is defined, for $x\in\mathbb{R}$ and $\nu\in\mathbb{R}$, by the power series
\begin{equation*}\mathbf{L}_\nu(x)=\sum_{k=0}^\infty \frac{\big(\frac{1}{2}x\big)^{\nu+2k+1}}{\Gamma(k+\frac{3}{2})\Gamma(k+\nu+\frac{3}{2})}.
\end{equation*}
The modified Struve function $\mathbf{L}_\nu(x)$ is closely related to the modified Bessel function $I_\nu(x)$, either sharing or having close analogues to the properties of $I_\nu(x)$ that were used by \cite{gaunt ineq1,gaunt ineq3,gaunt ineq8} to derive inequalities for the integral (\ref{intbes}).  The function $\mathbf{L}_\nu(x)$ is itself a widely used special function, with numerous applications in the applied sciences, such as perturbation approximations of lee waves in a stratified flow \cite{mh69}, leakage inductance in transformer windings \cite{hw94}, and quantum-statistical distribution functions of a hard-sphere system \cite{ni67}; see \cite{bp13} for examples of further application areas.  Basic properties of the modified Struve function  $\mathbf{L}_\nu(x)$ can be found in standard references,  such as \cite{olver}. We collect the basic properties that will be needed in this paper in Appendix \ref{appa}

The natural analogue of the problem studied by \cite{gaunt ineq1,gaunt ineq3,gaunt ineq8} is to ask for simple inequalities, involving the modified Struve function of the first kind, for the integral
\begin{equation}\label{intstruve}\int_0^x \mathrm{e}^{-\beta t} t^{\nu} \mathbf{L}_\nu(t)\,\mathrm{d}t, 
\end{equation}
where $x>0$, $0\leq\beta<1$ and $\nu>-1$ (with the condition on $\nu$ ensuring the integral exists).  This problem was first studied in the recent paper \cite{gaunt ineq4}, and will also be the subject of this paper.

The integral (\ref{intstruve}) can be evaluated exactly in terms of the modified Struve function $\mathbf{L}_\nu(x)$ in the case $\beta=1$. For all $\nu>-\frac{1}{2}$ and $x>0$,
\begin{equation}\label{intfor}\int_0^x \mathrm{e}^{-t}t^\nu \mathbf{L}_\nu(t)\,\mathrm{d}t=\frac{\mathrm{e}^{-x}x^{\nu+1}}{2\nu+1}\big(\mathbf{L}_\nu(x)+\mathbf{L}_{\nu+1}(x)\big)-\frac{\gamma(2\nu+2,x)}{\sqrt{\pi}2^\nu(2\nu+1)\Gamma(\nu+\frac{3}{2})},
\end{equation}
where  $\gamma(a,x)=\int_0^x \mathrm{e}^{-t}t^{a-1}\,\mathrm{d}t$ is the lower incomplete gamma function. 
This formula can be verified directly by a short calculation using the differentiation formula (\ref{diffone}) and identity (\ref{Iidentity}) given in Appendix \ref{appa}.
When $\beta=0$ the integral (\ref{intstruve}) cannot be evaluated in terms of the function $\mathbf{L}_\nu(x)$, but an exact formula is available in terms of the generalized hypergeometric function
\begin{equation*}{}_pF_q\big(a_1,\ldots,a_p;b_1,\ldots,b_q;x\big)=\sum_{k=0}^\infty \frac{(a_1)_k\cdots(a_p)_k}{(b_1)_k\cdots(b_q)_k}\frac{x^k}{k!},
\end{equation*}
where the Pochhammer symbol is defined by $(a)_0=1$ and $(a)_k=a(a+1)(a+2)\cdots(a+k-1)$, $k\geq1$.  Indeed, for $-\nu-\frac{3}{2}\notin\mathbb{N}$, we have the representation
\begin{equation*}\mathbf{L}_\nu(x)=\frac{x^{\nu+1}}{\sqrt{\pi}2^\nu\Gamma(\nu+\frac{3}{2})} {}_1F_2\bigg(1;\frac{3}{2},\nu+\frac{3}{2};\frac{x^2}{4}\bigg),
\end{equation*}
and by a straightforward calculation we have that, for $\nu>-1$ and $x>0$,
\begin{equation*}\int_0^x t^\nu \mathbf{L}_\nu(t)\,\mathrm{d}t=\frac{x^{2\nu+2}}{\sqrt{\pi}2^{\nu+1}(\nu+1)\Gamma(\nu+\frac{3}{2})}{}_2F_3\bigg(1,\nu+1;\frac{3}{2},\nu+\frac{3}{2},\nu+2;\frac{x^2}{4}\bigg).
\end{equation*}
   The integral (\ref{intstruve}) can also be evaluated when $0<\beta<1$, but the formula is more complicated: for $\nu>-1$ and $x>0$,
\begin{equation*}\int_0^x \mathrm{e}^{-\beta t} t^{\nu}  \mathbf{L}_{\nu}(t)\,\mathrm{d}t=\sum_{k=0}^\infty\frac{2^{-\nu-2k}\beta^{-2k-2\nu-2}}{\Gamma(k+\frac{3}{2})\Gamma(k+\nu+\frac{3}{2})}\gamma(2k+2\nu+2,\beta x).
\end{equation*}
  These complicated formulas provide the motivation for establishing simple bounds, involving the modified Struve function $\mathbf{L}_{\nu}(x)$ itself, for the integral (\ref{intstruve}).

Several upper bounds and a lower bound for the integral (\ref{intstruve}) were established by \cite{gaunt ineq4} by adapting the techniques used by \cite{gaunt ineq1,gaunt ineq3} to bound the analogous integral (\ref{intbes}) involving the modified Bessel function $I_\nu(x)$.  In this paper, we complement the work of \cite{gaunt ineq4} by obtaining several lower bounds for the integral (\ref{intstruve}) (Theorem \ref{tiger1}), one of which is a strict improvement on the only lower bound given in \cite{gaunt ineq4}.  In fact, all lower bounds obtained in this paper are tight in the limit $x\rightarrow\infty$, a feature not seen in the lower bound of \cite{gaunt ineq4}. We also extend the range of validity of the upper bounds given in \cite{gaunt ineq4} from $\nu\geq\frac{1}{2}$ to $\nu>-\frac{1}{2}$ (Theorem \ref{tiger2}), with our bounds taking the same functional form, but with larger numerical constants.  We shall proceed in a similar manner to \cite{gaunt ineq4}, by adapting the approach used in the recent paper \cite{gaunt ineq8} to obtain similar improvements on the bounds of \cite{gaunt ineq1, gaunt ineq3} that were obtained for the related integral (\ref{intbes}) involving $I_\nu(x)$. We establish our upper bounds by proving a series of lemmas, which may be of independent interest.  Lemma \ref{lem1} gives another upper bound for the integral (\ref{intstruve}), which outperforms our bounds from Theorem \ref{tiger2} for `large' values of $x$. In Lemma \ref{lem0}, we provide a new bound for the ratio $\mathbf{L}_{\nu}(x)/\mathbf{L}_{\nu-1}(x)$.  Lemma \ref{lem2} gives monotonicity results and inequalities for some products involving the modified Struve function $\mathbf{L}_\nu(x)$ and the modified Bessel function of the second kind $K_\nu(x)$ that complement existing results concerning products involving the modified Bessel functions $I_\nu(x)$ and $K_\nu(x)$. The lemmas are collected and proved in Section \ref{seclem}, and the main results are proved in Section \ref{sec3}.
 Elementary properties of the modified Struve function $\mathbf{L}_\nu(x)$ and the modified Bessel functions that are needed in the paper are collected in Appendix \ref{appa}.




\section{Main results and comparisons}\label{sec2}

The inequalities given in the following Theorems \ref{tiger1} and \ref{tiger2} are natural analogues of inequalities that have been recently obtained by \cite{gaunt ineq8} for the related integral $\int_0^x \mathrm{e}^{-\beta t}t^{\nu}I_\nu(t)\,\mathrm{d}t$.  The inequalities also complement and improve on bounds of \cite{gaunt ineq4} for the integral (\ref{intstruve}). Theorems \ref{tiger1} and \ref{tiger2} and Proposition \ref{propone} below are proved in Section \ref{sec3}.



\begin{theorem}\label{tiger1}Let $0<\beta<1$. Then, for $x>0$,
\begin{align}\int_0^x \mathrm{e}^{-\beta t}t^{\nu}\mathbf{L}_\nu(t)\,\mathrm{d}t&>\frac{1}{1-\beta}\bigg\{\mathrm{e}^{-\beta x}x^\nu\mathbf{L}_\nu(x)\nonumber\\
\label{ineqb2}&\quad
-\frac{\gamma(2\nu+1,\beta x)}{\sqrt{\pi}2^\nu\beta^{2\nu+1}\Gamma(\nu+\frac{3}{2})}\bigg\}, \quad-\tfrac{1}{2}<\nu\leq0, \\
\int_0^x \mathrm{e}^{-\beta t}t^{\nu}\mathbf{L}_\nu(t)\,\mathrm{d}t&>\frac{1}{1-\beta}\bigg\{\bigg(1-\frac{4\nu^2}{(2\nu-1)(1-\beta)}\frac{1}{x}\bigg)\mathrm{e}^{-\beta x}x^\nu \mathbf{L}_\nu(x)\nonumber\\
\label{ineqb3}&\quad-\frac{\gamma(2\nu+1,\beta x)}{\sqrt{\pi}2^\nu\beta^{2\nu+1}\Gamma(\nu+\frac{3}{2})}\bigg\}, \quad \nu\geq\tfrac{3}{2}, \\
\label{ineqb4}\int_0^x \mathrm{e}^{-\beta t}t^{\nu}\mathbf{L}_\nu(t)\,\mathrm{d}t&>\mathrm{e}^{-\beta x}x^\nu\sum_{k=0}^\infty \beta^k \mathbf{L}_{\nu+k+1}(x), \quad \nu>-1.
\end{align}
Inequalities (\ref{ineqb2})--(\ref{ineqb4}) are tight in the limit $x\rightarrow\infty$. Recall that $\gamma(a,x)=\int_0^x\mathrm{e}^{-t}t^{a-1}\,\mathrm{d}t$ is the lower incomplete gamma function.
\end{theorem}


\begin{theorem}\label{tiger2}Let $0<\beta<1$.  Then, for $x>0$,
\begin{align}\label{ineqb10}\int_0^x\mathrm{e}^{-\beta t}t^\nu \mathbf{L}_\nu(t)\,\mathrm{d}t&<\frac{2\nu+29}{(2\nu+1)(1-\beta)}\mathrm{e}^{-\beta x}x^\nu \mathbf{L}_{\nu+1}(x), \quad \nu>-\tfrac{1}{2}, \\
\label{ineqb11}\int_0^x\mathrm{e}^{-\beta t}t^\nu \mathbf{L}_\nu(t)\,\mathrm{d}t&<\frac{2\nu+15}{(2\nu+1)(1-\beta)}\mathrm{e}^{-\beta x}x^\nu \mathbf{L}_\nu(x), \quad \nu>-\tfrac{1}{2},
\end{align}
\begin{align}
\int_0^x\mathrm{e}^{-\beta t}t^\nu \mathbf{L}_\nu(t)\,\mathrm{d}t&>\frac{1}{1-\beta}\bigg\{\bigg(1-\frac{2\nu(2\nu+27)}{(2\nu-1)(1-\beta)}\frac{1}{x}\bigg)\mathrm{e}^{-\beta x}x^\nu \mathbf{L}_\nu(x)\nonumber\\
\label{ineqb12}&\quad-\frac{\gamma(2\nu+1,\beta x)}{\sqrt{\pi}2^\nu\beta^{2\nu+1}\Gamma(\nu+\frac{3}{2})}\bigg\}, \quad \nu>\tfrac{1}{2}.
\end{align}
Inequality (\ref{ineqb12}) is tight as $x\rightarrow\infty$. 
\end{theorem}

The inequalities in the following proposition are stronger than inequalities (\ref{ineqb2}), (\ref{ineqb3}) and (\ref{ineqb12}), because $\mathbf{L}_{\nu+1}(x)<\mathbf{L}_\nu(x)$, $x>0$, $\nu\geq-\frac{1}{2}$ (see (\ref{Imon})).

\begin{proposition}\label{propone}Let $0<\beta<1$. Then, for $x>0$,
\begin{align}\int_0^x \mathrm{e}^{-\beta t}t^{\nu}\mathbf{L}_{\nu+1}(t)\,\mathrm{d}t&>\frac{1}{1-\beta}\bigg\{\mathrm{e}^{-\beta x}x^\nu\mathbf{L}_\nu(x)\nonumber\\
\label{ineqb21}&\quad-\frac{\gamma(2\nu+1,\beta x)}{\sqrt{\pi}2^\nu\beta^{2\nu+1}\Gamma(\nu+\frac{3}{2})}\bigg\}, \quad-\tfrac{1}{2}<\nu\leq0,\\
\int_0^x \mathrm{e}^{-\beta t}t^{\nu}\mathbf{L}_{\nu+1}(t)\,\mathrm{d}t&>\frac{1}{1-\beta}\bigg\{\bigg(1-\frac{4\nu^2}{(2\nu-1)(1-\beta)}\frac{1}{x}\bigg)\mathrm{e}^{-\beta x}x^\nu \mathbf{L}_\nu(x)\nonumber\\
\label{ineqb22}&\quad-\frac{\gamma(2\nu+1,\beta x)}{\sqrt{\pi}2^\nu\beta^{2\nu+1}\Gamma(\nu+\frac{3}{2})}\bigg\}, \quad \nu\geq\tfrac{3}{2}, \\
\int_0^x\mathrm{e}^{-\beta t}t^\nu \mathbf{L}_{\nu+1}(t)\,\mathrm{d}t&>\frac{1}{1-\beta}\bigg\{\bigg(1-\frac{2\nu(2\nu+27)}{(2\nu-1)(1-\beta)}\frac{1}{x}\bigg)\mathrm{e}^{-\beta x}x^\nu \mathbf{L}_\nu(x)\nonumber\\
\label{ineqb23}&\quad -\frac{\gamma(2\nu+1,\beta x)}{\sqrt{\pi}2^\nu\beta^{2\nu+1}\Gamma(\nu+\frac{3}{2})}\bigg\}, \quad \nu>\tfrac{1}{2}.
\end{align}
\end{proposition}

\begin{remark}\label{rem1}In this remark, we discuss the performance of our bounds given in Theorems \ref{tiger1} and \ref{tiger2}, and make comparisons between our bounds and those given by \cite{gaunt ineq4} for the integral (\ref{intstruve}). Throughout this remark $0<\beta<1$. 

Inequality (\ref{ineqb4}) improves on the only other lower bound for the integral (\ref{intstruve}) in the literature \cite{gaunt ineq4}, $\int_0^x \mathrm{e}^{-\beta t}t^{\nu}\mathbf{L}_\nu(t)\,\mathrm{d}t>\mathrm{e}^{-\beta x}x^\nu \mathbf{L}_{\nu+1}(x)$, $x>0$, $\nu>-\tfrac{1}{2}$, with this bound in fact being the first term in the infinite series of the lower bound (\ref{ineqb4}). The other lower bounds from Theorems \ref{tiger1} and \ref{tiger2}, that is (\ref{ineqb2}), (\ref{ineqb3}) and (\ref{ineqb12}), all perform worse than (\ref{ineqb4}) and the bound of \cite{gaunt ineq4} for `small' $x$. Indeed, it is easily seen that the lower bounds in (\ref{ineqb2}) and (\ref{ineqb3}) are negative for sufficiently small $x$, whilst a simple asymptotic analysis of the bound (\ref{ineqb12}) using (\ref{Itend0}) shows that, for $-\frac{1}{2}<\nu<0$, the limiting form of this bound is $\frac{2\nu}{2\nu+1}\frac{x^{2\nu+1}}{\sqrt{\pi}2^\nu\Gamma(\nu+3/2)}<0$, as $x\downarrow0$. For the case $\nu=0$ the bound is again negative for sufficiently small $x$: $\frac{1}{1-\beta}\{\mathrm{e}^{-\beta x}\mathbf{L}_0(x)-\frac{2}{\pi\beta}(1-\mathrm{e}^{-\beta x})\}\sim -\frac{\beta x^2}{\pi(1-\beta)}$, as $x\downarrow0$.  The bounds (\ref{ineqb2}), (\ref{ineqb3}) and (\ref{ineqb12}) do, however, perform well for `large' $x$. Unlike the bound of \cite{gaunt ineq4}, these bounds are tight as $x\rightarrow\infty$, and this is achieved without the need of an infinite sum involving modified Struve functions of the first kind as given in the bound (\ref{ineqb4}).

 Inequality (2.13) of \cite{gaunt ineq4} gives the following upper bound: for $x>0$,
\begin{align}\int_0^x \mathrm{e}^{-\beta t}t^{\nu}\mathbf{L}_\nu(t)\,\mathrm{d}t&<\frac{\mathrm{e}^{-\beta x}x^\nu}{(2\nu+1)(1-\beta)}\bigg(2(\nu+1)\mathbf{L}_{\nu+1}(x)-\mathbf{L}_{\nu+3}(x) \nonumber\\
&\quad-\frac{x^{\nu+2}}{\sqrt{\pi}2^{\nu+2}(\nu+1)\Gamma(\nu+\frac{5}{2})}\bigg), \quad \nu\geq\tfrac{1}{2},\nonumber 
\end{align}
\begin{align}
\label{gau1}&<\frac{2(\nu+1)}{(2\nu+1)(1-\beta)}\mathrm{e}^{-\beta x}x^\nu \mathbf{L}_{\nu+1}(x), \quad \nu\geq\tfrac{1}{2}.
\end{align}
Another upper bound is obtained by combining inequalities (2.10) and (2.12) of \cite{gaunt ineq4}: for $x>0$,
\begin{equation}\label{gau2}\int_0^x \mathrm{e}^{-\beta t}t^{\nu}\mathbf{L}_\nu(t)\,\mathrm{d}t<\frac{1}{1-\beta}\mathrm{e}^{-\beta x}x^\nu \mathbf{L}_\nu(x), \quad \nu\geq\tfrac{1}{2}.
\end{equation}
Inequalities (\ref{ineqb10}) and (\ref{ineqb11}) increase the range of validity of inequalities (\ref{gau1}) and (\ref{gau2}) to $\nu>-\frac{1}{2}$ at the cost of larger multiplicative constants. These larger constants arise because our derivations of inequalities (\ref{ineqb10}) and (\ref{ineqb11}) are more involved than those of \cite{gaunt ineq4} for inequalities (\ref{gau1}) and (\ref{gau2}).  Indeed, we arrive at our bounds by applying a series of inequalities collected in Lemmas \ref{lem0}--\ref{notfin}, which when combined leads to a build up of errors. The reason we needed a more involved analysis was because the derivations of \cite{gaunt ineq4} rely heavily on the use of the inequality $\mathbf{L}_\nu(x)<\mathbf{L}_{\nu-1}(x)$, which holds for $x>0$, $\nu\geq\frac{1}{2}$ (see (\ref{Imon})), and without this useful inequality at our disposal (we have $\nu>-\frac{1}{2}$) we required a more involved and less direct proof. It is worth noting that we can combine our bound (\ref{ineqb10}) and the bound (\ref{gau1}) of \cite{gaunt ineq4} to obtain the bound, for $x>0$,
\begin{equation*}\int_0^x \mathrm{e}^{-\beta t}t^{\nu}\mathbf{L}_\nu(t)\,\mathrm{d}t<\frac{A_\nu}{(2\nu+1)(1-\beta)}\mathrm{e}^{-\beta x}x^\nu \mathbf{L}_{\nu+1}(x), \quad \nu>-\tfrac{1}{2},
\end{equation*}
where $A_\nu=2(\nu+1)$ for $\nu\geq\frac{1}{2}$, and $A_\nu=2\nu+29$ for $|\nu|<\frac{1}{2}$. A similar inequality can be obtained by combining our bound (\ref{ineqb11}) and the bound (\ref{gau2}) of \cite{gaunt ineq4}.


The inequalities obtained in this paper along with those presented in this remark allow for various double inequalities to be given for the integral (\ref{intstruve}).
As an example, for $x>0$,
\begin{equation}\label{gau3}\mathrm{e}^{-\beta x}x^\nu\sum_{k=0}^\infty \beta^k \mathbf{L}_{\nu+k+1}(x)<\int_0^x \mathrm{e}^{-\beta t}t^{\nu}\mathbf{L}_\nu(t)\,\mathrm{d}t<\frac{1}{1-\beta}\mathrm{e}^{-\beta x}x^\nu \mathbf{L}_\nu(x), \quad \nu\geq\tfrac{1}{2}.
\end{equation}
With the aid of \emph{Mathematica} we calculated the relative error in estimating  $F_{\nu,\beta}(x)=\int_0^x \mathrm{e}^{-\beta t}t^{\nu}\mathbf{L}_\nu(t)\,\mathrm{d}t$ by the upper bound in (\ref{gau3}) (denoted by $U_{\nu,\beta}(x)$), and the lower bound truncated at the fifth term in the sum,  $L_{\nu,\beta}(x)=\mathrm{e}^{-\beta x}x^\nu\sum_{k=0}^4 \beta^k \mathbf{L}_{\nu+k+1}(x)$. We report the results in Tables \ref{table1} and \ref{table2}. For fixed $x$ and $\nu$, we see that increasing $\beta$ increases the relative error in approximating $F_{\nu,\beta}(x)$ by either $L_{\nu,\beta}(x)$ or $U_{\nu,\beta}(x)$.  Both the lower and upper bounds in (\ref{gau3}) are tight as $x\rightarrow\infty$, and we see that, for fixed $\nu$ and $\beta$, the relative error in approximating $F_{\nu,\beta}(x)$ by $U_{\nu,\beta}(x)$ decreases as $x$ increases. However, as we have truncated the sum, $L_{\nu,\beta}(x)$ is not tight as $x\rightarrow\infty$.  The effect of truncating the sum is most pronounced for larger $\beta$ and larger $x$. For $\beta=0.75$, $\sum_{k=0}^\infty 0.75^k=4$ and $\sum_{k=0}^4 0.75^k=3.0508$, and so $\lim_{x\rightarrow\infty}\big(1-\frac{L_{\nu,0.75}(x)}{F_{\nu,0.75}(x)}\big)=0.2373$, $\nu>-\frac{1}{2}$, where we also made use of the limiting forms (\ref{eqeq1}) and (\ref{Itendinfinity}). In contrast, $\lim_{x\rightarrow\infty}\big(1-\frac{L_{\nu,0.25}(x)}{F_{\nu,0.25}(x)}\big)=9.766\times 10^{-4}$, which is fairly negligible. The upper bound $U_{\nu,\beta}(x)$ is of the wrong asymptotic order as $x\downarrow0$ (using (\ref{Itend0}) shows that $\frac{U_{\nu,\beta}(x)}{F_{\nu,\beta}(x)}\sim\frac{2(\nu+1)}{(1-\beta)x}$, as $x\downarrow0$), and so performs poorly for `small' $x$.  The lower bound $L_{\nu,\beta}(x)$ performs better for `small' $x$; indeed, it is of the correct asymptotic order as $x\downarrow0$ with $\lim_{x\downarrow0}\big(1-\frac{L_{\nu,\beta}(x)}{F_{\nu,\beta}(x)}\big)=\frac{1}{2\nu+3}$.  

\begin{table}[h]
\centering
\caption{\footnotesize{Relative error in approximating $F_{\nu,\beta}(x)$ by $L_{\nu,\beta}(x)$.}}
\label{table1}
{\scriptsize
\begin{tabular}{|c|rrrrrrr|}
\hline
 \backslashbox{$(\nu,\beta)$}{$x$}      &    0.5 &    5 &    10 &    15 &    25 &    50 & 100   \\
 \hline
$(1,0.25)$ & 0.2051 & 0.1976 & 0.1413 & 0.1028 & 0.0656 & 0.0346 & 0.0182 \\
$(2.5,0.25)$ & 0.1276 & 0.1320 & 0.1092 & 0.0863 & 0.0591 & 0.0329 &  0.0177 \\
$(5,0.25)$ & 0.0781 & 0.0831 & 0.0773 & 0.0670 & 0.0503 & 0.0302 & 0.0169  \\
$(10,0.25)$ & 0.0439 & 0.0465 & 0.0468 & 0.0444 & 0.0378  & 0.0257 & 0.0155  \\   
  \hline
$(1,0.5)$ & 0.2111 & 0.2582 & 0.2259 & 0.1843 & 0.1341 & 0.0870 & 0.0602 \\
$(2.5,0.5)$ & 0.1304 & 0.1635 & 0.1606 & 0.1426 & 0.1133 & 0.0791 &   0.0570 \\
$(5,0.5)$ & 0.0793 & 0.0971 & 0.1039 & 0.1004 & 0.0881 & 0.0680 & 0.0522 \\
$(10,0.5)$ & 0.0443 & 0.0514 & 0.0569  & 0.0590 &  0.0580 & 0.0515 & 0.0440  \\ 
\hline
$(1,0.75)$ & 0.2171 & 0.3359 & 0.3723 & 0.3659 & 0.3369 & 0.2953 &  0.2683  \\
$(2.5,0.75)$ & 0.1333 & 0.2036 & 0.2458 & 0.2597 & 0.2640 & 0.2581 &   0.2500 \\
$(5,0.75)$ & 0.0805 & 0.1142 & 0.1446 & 0.1635 & 0.1850 & 0.2084 & 0.2226  \\ 
$(10,0.75)$ & 0.0447 & 0.0569 & 0.0705 & 0.0825 & 0.1028 & 0.1400 & 0.1774 \\  
  \hline
\end{tabular}}
\end{table}

\begin{table}[h]
\centering
\caption{\footnotesize{Relative error in approximating $F_{\nu,\beta}(x)$ by $U_{\nu,\beta}(x)$.}}
\label{table2}
{\scriptsize
\begin{tabular}{|c|rrrrrrr|}
\hline
 \backslashbox{$(\nu,\beta)$}{$x$}      &    0.5 &    5 &    10 &    15 &    25 &    50 & 100   \\
 \hline
$(1,0.25)$ & 9.4597 & 0.3208 & 0.0888 & 0.0521 & 0.0292 & 0.0139 & 0.0068 \\
$(2.5,0.25)$ & 17.4185 & 0.9887 & 0.3593 & 0.2156 & 0.1197  & 0.0565 & 0.0274 \\
$(5,0.25)$ & 30.7218 & 2.1879 & 0.8593 & 0.5134 & 0.2806 & 0.1300 & 0.0625  \\
$(10,0.25)$ & 57.3655 & 4.7301 & 1.9918 & 1.1901 & 0.6378  & 0.2868 & 0.1351 \\   
  \hline
$(1,0.5)$ & 14.2938 & 0.5538 & 0.1530 & 0.0839 & 0.0452 & 0.0212 & 0.0103 \\
$(2.5,0.5)$ & 26.1923 & 1.5400 & 0.5661 & 0.3363 & 0.1842 & 0.0858 &  0.0414  \\
$(5,0.5)$ & 46.1220 & 3.3214 & 1.3161 & 0.7868 & 0.4286 & 0.1972 & 0.0943 \\
$(10,0.5)$ & 86.0701 & 7.1185 & 3.0084  & 1.8015 & 0.9664  & 0.4339 & 0.2037 \\ 
\hline
$(1,0.75)$ & 28.8028 & 1.3243 & 0.4124 & 0.2137 & 0.1021 & 0.0444 &   0.0210 \\
$(2.5,0.75)$ & 52.5169 & 3.2293 & 1.2300 & 0.7305 & 0.3933 & 0.1783 &  0.0845 \\
$(5,0.75)$ & 92.3236 & 6.7374 & 2.7112 & 1.6308 & 0.8892 & 0.4056 & 0.1918  \\
$(10,0.75)$ & 172.1854 & 14.2887 & 6.0686 & 3.6482 & 1.9648  & 0.8827 & 0.4126 \\  
  \hline
\end{tabular}}
\end{table}
\end{remark}

\section{Lemmas}\label{seclem}


We prove Theorem \ref{tiger2} through the following series of lemmas, which may be of independent interest. 

\begin{lemma}\label{lem0}Let $\nu>0$ and $x>0$. Then
\begin{align}\label{lratio}\frac{\mathbf{L}_{\nu}(x)}{\mathbf{L}_{\nu-1}(x)}>\frac{x}{2\nu+1+x}.
\end{align}
This bound is tight in the limits $x\downarrow0$ and $x\rightarrow\infty$.
\end{lemma}

\begin{lemma}\label{lem2}Suppose $\nu\geq-\frac{1}{2}$. Then the functions $x\mapsto K_{\nu+1}(x)\mathbf{L}_\nu(x)$ and $x\mapsto xK_{\nu+2}(x)\mathbf{L}_\nu(x)$ are strictly decreasing on $(0,\infty)$. As a consequence of the latter monotonicity result, we have the following tight two-sided inequality: 
\begin{equation}\label{klineq1}\frac{1}{2}<xK_{\nu+2}(x)\mathbf{L}_\nu(x)<\frac{2\Gamma(\nu+2)}{\sqrt{\pi}\Gamma(\nu+\frac{3}{2})},\quad x>0.
\end{equation}
We also have that, for $x>0$,
\begin{align}\label{klineq0}xK_{\nu+1}(x)\mathbf{L}_\nu(x)&<1,\\
\label{klineq2}xK_{\nu+3}(x)\mathbf{L}_\nu(x)&<\frac{2\Gamma(\nu+2)}{\sqrt{\pi}\Gamma(\nu+\frac{3}{2})}\bigg(1+\frac{2\nu+5}{x}\bigg).
\end{align}
Suppose now that $-\frac{1}{2}\leq\nu\leq\frac{1}{2}$. Then, for $x>0$,
\begin{align}\label{gamr1}xK_{\nu+2}(x)\mathbf{L}_\nu(x)&<\frac{3}{2},\\
\label{gamr2}xK_{\nu+3}(x)\mathbf{L}_\nu(x)&<\frac{3}{2}+\frac{9}{x},\\
\label{gamr3}xK_{\nu+3}(x)\mathbf{L}_{\nu+1}(x)&<\frac{15}{8}.
\end{align}
\end{lemma}

\begin{lemma}\label{lem1}Let $\nu>-\frac{1}{2}$ and $0<\beta<1$. Fix $x_*>\frac{1}{1-\beta}$.  Then, for $x\geq x_*$,
\begin{equation}\label{ineqb1}\int_0^x \mathrm{e}^{-\beta t}t^{\nu}\mathbf{L}_\nu(t)\,\mathrm{d}t<M_{\nu,\beta}(x_*)\mathrm{e}^{-\beta x} x^\nu \mathbf{L}_{\nu+1}(x), 
\end{equation}
where
\begin{equation}\label{mng}M_{\nu,\beta}(x_*)=\max\bigg\{\frac{2\nu+3+2x_*}{2\nu+1},\frac{x_*}{(1-\beta)x_*-1}\bigg\}.
\end{equation}
\end{lemma}

\begin{lemma}\label{notfin} Suppose that $-\frac{1}{2}<\nu\leq\frac{1}{2}$ and $0<\beta<1$. Then, for $x>0$,
\begin{align}\label{term00}\frac{\mathrm{e}^{\beta x}K_{\nu+3}(x)}{x^{\nu-1}}\int_0^x \mathrm{e}^{-\beta t}t^{\nu}\mathbf{L}_\nu(t)\,\mathrm{d}t&<\frac{14}{(2\nu+1)(1-\beta)}, \\
\label{term10}\frac{\mathrm{e}^{\beta x}K_{\nu+2}(x)}{x^{\nu-1}}\int_0^x \mathrm{e}^{-\beta t}t^{\nu}\mathbf{L}_\nu(t)\,\mathrm{d}t&<\frac{7}{(2\nu+1)(1-\beta)}.
\end{align}
\end{lemma}


\begin{remark}The monotonicity results of Lemma \ref{lem2} for the products $K_{\nu+1}(x)\mathbf{L}_\nu(x)$ and $xK_{\nu+2}(x)\mathbf{L}_\nu(x)$ complement monotonicity results that have been established for the products $K_\nu(x)I_\nu(x)$ (see \cite{bar1,bar2,penfold,pm50}), $xK_\nu(x)I_\nu(x)$ (see \cite{hartman}) and $xK_{\nu+1}(x)I_\nu(x)$ (see \cite{gaunt ineq2}).  We also note that a number of bounds for the product $K_\nu(x)I_\nu(x)$ have been obtained by \cite{bar16}.  In light of these results, it is natural to ask whether a monotonicity result is available for the product $xK_{\nu+1}(x)\mathbf{L}_\nu(x)$, which is also present in Lemma \ref{lem2}.  It turns out that, for fixed $\nu>-\frac{1}{2}$, $xK_{\nu+1}(x)\mathbf{L}_\nu(x)$ is not a monotone function of $x$ on $(0,\infty)$. Indeed, applying the limiting forms (\ref{Itend0})--(\ref{Ktendinfinity}) gives that
\begin{align*}xK_{\nu+1}(x)\mathbf{L}_\nu(x) &\sim\frac{\Gamma(\nu+1)x}{\sqrt{\pi}\Gamma(\nu+\frac{3}{2})}, \quad x\downarrow0, \\
xK_{\nu+1}(x)\mathbf{L}_\nu(x) &\sim \frac{1}{2}+\frac{2\nu+1}{4x}, \quad x\rightarrow\infty,
\end{align*}
which tells us that $xK_{\nu+1}(x)\mathbf{L}_\nu(x)$ is an increasing function of $x$ for `small' $x$ and a decreasing function of $x$ for `large' $x$ if $\nu>-\frac{1}{2}$.
\end{remark}

\begin{remark}Inequality (\ref{ineqb1}) of Lemma \ref{lem1} is more accurate than inequalities (\ref{ineqb10}) and (\ref{ineqb11}) of Theorem \ref{tiger1} for `large' $x$.  As an example, applying Lemma \ref{lem1} with $x_*=\frac{2}{1-\beta}$ gives that, for $x\geq\frac{2}{1-\beta}$, $\nu>-\frac{1}{2}$, $0<\beta<1$,
\begin{equation*}\int_0^x\mathrm{e}^{-\beta t}t^\nu\mathbf{L}_\nu(t)\,\mathrm{d}t<\frac{1}{2\nu+1}\bigg(2\nu+3+\frac{4}{1-\beta}\bigg)\mathrm{e}^{-\beta x}x^\nu\mathbf{L}_{\nu+1}(x),
\end{equation*} 
which is an improvement on both (\ref{ineqb10}) and (\ref{ineqb11}) in its range of validity.
\end{remark}

\noindent{\emph{Proof of Lemma \ref{lem0}.}} We begin by noting the following bound of \cite[Theorem 2.2]{gaunt ineq5}:  
\begin{equation*}\frac{\mathbf{L}_\nu(x)}{\mathbf{L}_{\nu-1}(x)}>\bigg(\frac{I_{\nu-1}(x)}{I_\nu(x)}+\frac{2b_\nu(x)}{x}\bigg)^{-1}, \quad x>0,\:\nu\geq0,
\end{equation*}
where $b_\nu(x)=\frac{(x/2)^{\nu+1}}{\sqrt{\pi}\Gamma(\nu+3/2)\mathbf{L}_\nu(x)}$. 
Part (iii) of Lemma 2.1 of \cite{gaunt ineq5} tells us that $b_\nu(x)<\frac{1}{2}$ for all $x>0$, $\nu\geq0$, and so we have the simpler bound
\begin{equation}\label{aug18}\frac{\mathbf{L}_\nu(x)}{\mathbf{L}_{\nu-1}(x)}>\bigg(\frac{I_{\nu-1}(x)}{I_\nu(x)}+\frac{1}{x}\bigg)^{-1}, \quad x>0,\;\nu\geq0.
\end{equation}
The ratio of modified Bessel functions of the first kind can be bounded by the inequality 
\begin{equation*}\frac{I_{\nu}(x)}{I_{\nu-1}(x)}>\frac{x}{2\nu+x}, \quad x>0,\;\nu>0,
\end{equation*}
which is the simplest lower bound in a sequence of rational bounds obtained by \cite{nasell2}. Applying this bound to (\ref{aug18}) then gives us our desired bound (\ref{lratio}). Finally, the assertion that the bound is tight in the limits $x\downarrow0$ and $x\rightarrow\infty$ follow easily from an application of the limiting forms (\ref{Itend0}) and (\ref{Itendinfinity}) and the standard formula $\Gamma(x+1)=x\Gamma(x)$. \hfill $\square$

\vspace{3mm}

\noindent{\emph{Proof of Lemma \ref{lem2}.}} (i) Note that we can write $K_{\nu+1}(x)\mathbf{L}_\nu(x)=f_\nu(x)g_\nu(x)$, where $f_\nu(x)=K_{\nu+1}(x)I_{\nu+1}(x)$ and $g_\nu(x)=\mathbf{L}_\nu(x)/I_{\nu+1}(x)$. It has been shown that, for $\nu>-2$, $f_\nu(x)$ is a strictly decreasing function of $x$ on $(0,\infty)$ (see \cite{bar2}, which extends the range of validity of results of \cite{bar1,penfold}), and part (i) of Theorem 2.2 of \cite{bp14} states that, for $\nu\geq-\frac{1}{2}$, $g_\nu(x)$ is a decreasing function of $x$ on $(0,\infty)$. As a product of two strictly positive functions, one of which is strictly decreasing and the other decreasing, it follows that, for $\nu\geq-\frac{1}{2}$, the function $x\mapsto K_{\nu+1}(x)\mathbf{L}_\nu(x)$ is strictly decreasing on $(0,\infty)$. 

The proof that, for $\nu\geq-\frac{1}{2}$, the function $x\mapsto xK_{\nu+2}(x)\mathbf{L}_\nu(x)$ is strictly decreasing on $(0,\infty)$ is similar.  We note that $xK_{\nu+2}(x)\mathbf{L}_\nu(x)=h_\nu(x)g_\nu(x)$, where $h_\nu(x)=xK_{\nu+2}(x)I_{\nu+1}(x)$. Lemma 3 of \cite{gaunt ineq2} asserts that, for $\nu\geq-\frac{3}{2}$, $f_\nu(x)$ is a strictly decreasing function of $x$ on $(0,\infty)$, and the proof now proceeds exactly as the previous one concerning the monotonicity of the function $x\mapsto K_{\nu+1}(x)\mathbf{L}_\nu(x)$. The upper and lower bounds in (\ref{klineq1}) now follow from using the limiting forms (\ref{Itend0})--(\ref{Ktendinfinity}) to calculate the limits $\lim_{x\downarrow0}xK_{\nu+2}(x)\mathbf{L}_\nu(x)$ and $\lim_{x\rightarrow\infty}xK_{\nu+2}(x)\mathbf{L}_\nu(x)$. 


\vspace{1mm}

\noindent{(ii)} Inequality (\ref{klineq0}) is obtained by combining the inequality $\mathbf{L}_\nu(x)<I_\nu(x)$, $x>0$, $\nu\geq-\frac{1}{2}$, with the bound $xK_{\nu+1}(x)I_\nu(x)\leq1$, $x>0$, $\nu\geq-\frac{1}{2}$ (see \cite[Lemma 3]{gaunt ineq2}).  To see that $\mathbf{L}_\nu(x)<I_\nu(x)$, $x>0$, $\nu\geq-\frac{1}{2}$, we recall that the modified Struve function of the second kind is defined by $\mathbf{M}_\nu(x)=\mathbf{L}_\nu(x)-I_\nu(x)$. We can readily see that $\mathbf{M}_\nu(x)<0$, for $x>0$, $\nu>-\frac{1}{2}$, from its integral representation (see \cite[formula 11.5.4]{olver}), and $M_{-\frac{1}{2}}(x)<0$, $x>0$, can be seen by using the formulas in (\ref{speccase}). 

\vspace{1mm}

\noindent{(iii)} We will make use of the following inequality of \cite{segura} for a ratio of modified Bessel functions of the second kind:
\begin{equation}\label{kmuineq}\frac{K_\nu(x)}{K_{\nu-1}(x)}<\frac{\nu-\frac{1}{2}+\sqrt{(\nu-\frac{1}{2})^2+x^2}}{x}<1+\frac{2\nu-1}{x}, \quad x>0,\: \nu>\tfrac{1}{2}.
\end{equation}
We now obtain inequality (\ref{klineq2}) by applying inequality (\ref{kmuineq}) and the upper bound in (\ref{klineq1}):
\begin{align*}xK_{\nu+3}(x)\mathbf{L}_\nu(x)=\frac{K_{\nu+3}(x)}{K_{\nu+2}(x)}\cdot xK_{\nu+2}(x)\mathbf{L}_\nu(x)<\bigg(1+\frac{2\nu+5}{x}\bigg)\frac{2\Gamma(\nu+2)}{\sqrt{\pi}\Gamma(\nu+\frac{3}{2})}.
\end{align*}

\vspace{1mm}

\noindent{(iv)} We note that the ratio $\frac{\Gamma(\nu+2)}{\Gamma(\nu+3/2)}$ is an increasing function of $\nu$ on $[-\frac{1}{2},\frac{1}{2}]$ (see \cite{gio}). Therefore using the upper bound in (\ref{klineq1}) we obtain that, for $-\frac{1}{2}\leq\nu\leq\frac{1}{2}$ and $x>0$,
\[xK_{\nu+2}(x)\mathbf{L}_\nu(x)<\frac{2\Gamma(\frac{1}{2}+2)}{\sqrt{\pi}\Gamma(\frac{1}{2}+\frac{3}{2})}=\frac{3}{2},\]
where we used that $\Gamma(\frac{5}{2})=\frac{3\sqrt{\pi}}{4}$. Thus, we have proved inequality (\ref{gamr1}).  Inequalities (\ref{gamr2}) and (\ref{gamr3}) are obtained similarly (making use of the upper bound in (\ref{klineq1}) and inequality (\ref{klineq2})), and we omit the details.  \hfill $\square$

\vspace{3mm}

\noindent{\emph{Proof of Lemma \ref{lem1}.}} Fix $x_*>\frac{1}{1-\beta}$.
We consider the function
\begin{equation*}u_{\nu,\beta}(x)=M_{\nu,\beta}(x_*)\mathrm{e}^{-\beta x}x^\nu \mathbf{L}_{\nu+1}(x)-\int_0^x \mathrm{e}^{-\beta t}t^\nu \mathbf{L}_\nu(t)\,\mathrm{d}t,
\end{equation*}
and prove inequality (\ref{ineqb1}) by showing that $u_{\nu,\beta}(x)>0$ for all $x\geq x_*$.  

We first prove that $u_{\nu,\beta}(x_*)>0$. To this end, we consider the function
\begin{equation*}v_{\nu,\beta}(x)=\frac{\mathrm{e}^{\beta x}}{x^\nu \mathbf{L}_{\nu+1}(x)}\int_0^x\mathrm{e}^{-\beta t} t^\nu \mathbf{L}_\nu(t)\,\mathrm{d}t,
\end{equation*}
and it suffices to prove that $v_{\nu,\beta}(x_*)< M_{\nu,\beta}(x_*)$. We note that
\begin{equation*}\frac{\partial v_{\nu,\beta}(x)}{\partial \beta}=\frac{\mathrm{e}^{\beta x}}{x^\nu \mathbf{L}_{\nu+1}(x)}\int_0^x(x-t)\mathrm{e}^{-\beta t} t^\nu \mathbf{L}_\nu(t)\,\mathrm{d}t>0,
\end{equation*}
meaning that $v_{\nu,\beta}(x)$ is an increasing function of $\beta$. Therefore, for $0<\beta<1$,
\begin{align*}v_{\nu,\beta}(x_*)&< \frac{\mathrm{e}^{ x_*}}{x_*^\nu \mathbf{L}_{\nu+1}(x_*)}\int_0^{x_*}\mathrm{e}^{- t} t^\nu \mathbf{L}_\nu(t)\,\mathrm{d}t<\frac{x_*}{2\nu+1}\bigg(\frac{\mathbf{L}_\nu(x_*)}{\mathbf{L}_{\nu+1}(x_*)}+1\bigg) \\
&<\frac{x_*}{2\nu+1}\bigg(\frac{2\nu+3+x_*}{x_*}+1\bigg)=\frac{2\nu+3+2x_*}{2\nu+1}\leq M_{\nu,\beta}(x_*),
\end{align*}
where the second inequality is clear from the integral formula (\ref{intfor}) and we applied Lemma \ref{lem0} to obtain the third inequality. 

We now prove that $u_{\nu,\beta}'(x)>0$ for $x>x_*$.  A calculation using the differentiation formula (\ref{diffone}) followed  by an application of inequality (\ref{Imon}) gives that
\begin{align*}u_{\nu,\beta}'(x)&=M_{\nu,\beta}(x_*)\frac{\mathrm{d}}{\mathrm{d}x}\big(\mathrm{e}^{-\beta x}x^{-1}\cdot x^{\nu+1} \mathbf{L}_{\nu+1}(x)\big)- \mathrm{e}^{-\beta x}x^\nu \mathbf{L}_\nu(x)\\
&=M_{\nu,\beta}(x_*)\mathrm{e}^{-\beta x}x^\nu\big(\mathbf{L}_{\nu}(x)-x^{-1}\mathbf{L}_{\nu+1}(x)-\beta \mathbf{L}_{\nu+1}(x)\big)-\mathrm{e}^{-\beta x}x^\nu \mathbf{L}_\nu(x)\\
&>M_{\nu,\beta}(x_*)\mathrm{e}^{-\beta x}x^\nu\big(1-\beta -x^{-1})\mathbf{L}_\nu(x)-\mathrm{e}^{-\beta x}x^\nu \mathbf{L}_\nu(x)\\
&\geq\bigg(\frac{1-\beta-x^{-1}}{1-\beta-x_*^{-1}}-1\bigg)\mathrm{e}^{-\beta x}x^\nu \mathbf{L}_\nu(x)>0,
\end{align*}
for $x>x_*$. This completes the proof.
\hfill $\square$

\vspace{3mm}

\noindent{\emph{Proof of Lemma \ref{notfin}.}} (i) We obtain inequality (\ref{term00}) by bounding the expression
\[\frac{\mathrm{e}^{\beta x}K_{\nu+3}(x)}{x^{\nu-1}}\int_0^x \mathrm{e}^{-\beta t}t^{\nu}\mathbf{L}_\nu(t)\,\mathrm{d}t\]
for $x\in(0,x_*)$ and $x\in[x_*,\infty)$, where $x_*=\frac{C}{1-\beta}$ for some $C>1$ that we will choose later. Suppose first that $x\in(0,x_*)$. Observe that
\begin{equation*}\frac{\partial}{\partial \beta}\bigg(\frac{\mathrm{e}^{\beta x}K_{\nu+3}(x)}{x^{\nu-1}}\int_0^x \mathrm{e}^{-\beta t}t^{\nu}\mathbf{L}_\nu(t)\,\mathrm{d}t\bigg)=\frac{\mathrm{e}^{\beta x}K_{\nu+3}(x)}{x^{\nu-1}}\int_0^x (x-t)\mathrm{e}^{-\beta t}t^{\nu}\mathbf{L}_\nu(t)\,\mathrm{d}t> 0.
\end{equation*}
Since $0<\beta<1$, we therefore have that, for $x\in(0,x_*)$,
\begin{align}\frac{\mathrm{e}^{\beta x}K_{\nu+3}(x)}{x^{\nu-1}}\int_0^x \mathrm{e}^{-\beta t}t^{\nu}\mathbf{L}_\nu(t)\,\mathrm{d}t&< \frac{\mathrm{e}^{ x}K_{\nu+3}(x)}{x^{\nu-1}}\int_0^x \mathrm{e}^{- t}t^{\nu}\mathbf{L}_\nu(t)\,\mathrm{d}t\nonumber\\
&<\frac{1}{2\nu+1}x^2K_{\nu+3}(x)\big(\mathbf{L}_\nu(x)+\mathbf{L}_{\nu+1}(x)\big)\nonumber\\
&<\frac{1}{2\nu+1}\bigg(\frac{27}{8}x_*+9\bigg)=\frac{1}{2\nu+1}\bigg(9+\frac{27C}{8(1-\beta)}\bigg)\nonumber\\
&<\frac{1}{(2\nu+1)(1-\beta)}\bigg(9+\frac{27}{8}C\bigg)=:T_1,\nonumber
\end{align}
where we used (\ref{intfor}) to bound the integral in the second step, and inequalities (\ref{gamr2}) and (\ref{gamr3}) to obtain the third inequality.

Suppose now that $x\in[x_*,\infty)$. Let $M_{\nu,\beta}(x_*)$ be defined as per (\ref{mng}). Bounding the integral by inequality (\ref{ineqb1}) gives that
\begin{align}\frac{\mathrm{e}^{\beta x}K_{\nu+3}(x)}{x^{\nu-1}}\int_0^x \mathrm{e}^{-\beta t}&t^{\nu}\mathbf{L}_\nu(t)\,\mathrm{d}t<\frac{\mathrm{e}^{\beta x}K_{\nu+3}(x)}{x^{\nu-1}}\cdot M_{\nu,\beta}(x_*)\mathrm{e}^{-\beta x}x^\nu \mathbf{L}_{\nu+1}(x)\nonumber \\
&=M_{\nu,\beta}(x_*)xK_{\nu+3}(x)\mathbf{L}_{\nu+1}(x)\nonumber\\
&<\frac{15}{8} M_{\nu,\beta}(x_*)\nonumber\\
&=\frac{15}{8}\max\bigg\{\frac{1}{2\nu+1}\bigg(2\nu+3+\frac{2C}{1-\beta}\bigg),\frac{C}{(C-1)(1-\beta)}\bigg\}\nonumber\\
&\leq\max\bigg\{\frac{15(4+2C)}{8(2\nu+1)(1-\beta)},\frac{15C}{4(C-1)(2\nu+1)(1-\beta)}\bigg\}\nonumber\\
&=:\max\{T_2,T_3\},\nonumber 
\end{align}
where we used inequality (\ref{gamr3}) to obtain the second inequality and we used that $-\frac{1}{2}<\nu\leq\frac{1}{2}$ to obtain the third inequality.

It is readily checked that $T_1\geq T_2$ if $C\leq4$. Equating $T_1=T_3$ gives a quadratic equation for $C$ with positive solution $C=\frac{\sqrt{889}}{18}-\frac{5}{18}=1.3786\ldots$. Therefore
\begin{align*}\frac{\mathrm{e}^{\beta x}K_{\nu+3}(x)}{x^{\nu-1}}\int_0^x \mathrm{e}^{-\beta t}t^{\nu}\mathbf{L}_\nu(t)\,\mathrm{d}t&<\frac{1}{(2\nu+1)(1-\beta)}\bigg(9+\frac{27}{8}\cdot 1.3786\bigg)\\
&=\frac{13.653}{(2\nu+1)(1-\beta)}<\frac{14}{(2\nu+1)(1-\beta)}.
\end{align*}

\vspace{1mm}

\noindent{(ii)} The proof of inequality (\ref{term10}) is similar to that of inequality (\ref{term00}). Let $x_*=\frac{3}{2(1-\beta)}$. 
By a similar argument, we have that, for $x\in(0,x_*)$,
\begin{align}\frac{\mathrm{e}^{\beta x}K_{\nu+2}(x)}{x^{\nu-1}}\int_0^x \mathrm{e}^{-\beta t}t^{\nu}\mathbf{L}_\nu(t)\,\mathrm{d}t&<\frac{1}{2\nu+1}x^2K_{\nu+2}(x)\big(\mathbf{L}_\nu(x)+\mathbf{L}_{\nu+1}(x)\big)\nonumber\\
\label{lab1}&<\frac{x_*}{2\nu+1}\bigg(\frac{3}{2}+1\bigg)=\frac{15}{4(2\nu+1)(1-\beta)},
\end{align}
where we applied inequalities (\ref{gamr1}) and (\ref{klineq0}) to get the second inequality.
Suppose now that $x\in[x_*,\infty)$. Using inequality (\ref{ineqb1}) gives that
\begin{align}\frac{\mathrm{e}^{\beta x}K_{\nu+2}(x)}{x^{\nu-1}}\int_0^x \mathrm{e}^{-\beta t}t^{\nu}\mathbf{L}_\nu(t)\,\mathrm{d}t&<M_{\nu,\beta}(x_*)xK_{\nu+2}(x)\mathbf{L}_{\nu+1}(x)< M_{\nu,\beta}(x_*)\nonumber \\
&=\max\bigg\{\frac{1}{2\nu+1}\bigg(2\nu+3+\frac{3}{1-\beta}\bigg),\frac{3}{1-\beta}\bigg\}\nonumber\\
&\leq\max\bigg\{\frac{2\nu+6}{(2\nu+1)(1-\beta)},\frac{3}{1-\beta}\bigg\}\nonumber\\
\label{lab2}&=\frac{2\nu+6}{(2\nu+1)(1-\beta)}<\frac{7}{(2\nu+1)(1-\beta)},
\end{align}
where we used (\ref{klineq0}) to get the second inequality and we used that $-\frac{1}{2}<\nu\leq\frac{1}{2}$ to obtain the third and fourth inequalities.
We complete the proof by noting that the bound (\ref{lab2}) is greater than the bound (\ref{lab1}). \hfill $\square$

\section{Proofs of main results}\label{sec3}

\noindent{\emph{Proof of Theorem \ref{tiger1}.}} (i) Let $x>0$ and suppose $-\frac{1}{2}<\nu\leq0$.  Using integration by parts and the differentiation formula (\ref{diffone}) gives that
\begin{align*}\int_0^x\mathrm{e}^{-\beta t}t^\nu \mathbf{L}_\nu(t)\,\mathrm{d}t&=-\frac{1}{\beta}\mathrm{e}^{-\beta x}x^\nu \mathbf{L}_\nu(x)+\frac{1}{\beta}\int_0^x \mathrm{e}^{-\beta t}t^\nu \mathbf{L}_{\nu-1}(t)\,\mathrm{d}t,
\end{align*}
where we used that $\lim_{x\downarrow0}x^\nu\mathbf{L}_\nu(x)=0$, for $\nu>-\frac{1}{2}$ (see \ref{Itend0})). One can check that the integrals exist for $\nu>-\frac{1}{2}$ by using the limiting form (\ref{Itend0}).  By using the identity (\ref{Iidentity}) and rearranging we obtain that
\begin{align}&\int_0^x \mathrm{e}^{-\beta t}t^\nu \mathbf{L}_{\nu+1}(t)\,\mathrm{d}t+2\nu\int_0^x\mathrm{e}^{-\beta t}t^{\nu-1}\mathbf{L}_\nu(t)\,\mathrm{d}t-\beta\int_0^x\mathrm{e}^{-\beta t}t^\nu \mathbf{L}_\nu(t)\,\mathrm{d}t\nonumber\\
\label{55555}&\quad=\mathrm{e}^{-\beta x}x^\nu \mathbf{L}_\nu(x)-\int_0^x\mathrm{e}^{-\beta t}\frac{t^{2\nu}}{\sqrt{\pi}2^\nu\Gamma(\nu+\frac{3}{2})}\,\mathrm{d}t.
\end{align}
Using inequality (\ref{Imon}) to bound the first integral and making use of the assumption that $\nu\leq0$ gives that
\begin{align*}\int_0^x \mathrm{e}^{-\beta t}t^\nu \mathbf{L}_{\nu}(t)\,\mathrm{d}t&>\frac{1}{1-\beta}\bigg\{\mathrm{e}^{-\beta x}x^\nu \mathbf{L}_\nu(x)-\int_0^x\mathrm{e}^{-\beta t}\frac{t^{2\nu}}{\sqrt{\pi}2^\nu\Gamma(\nu+\frac{3}{2})}\,\mathrm{d}t\bigg\}.
\end{align*}
Finally, we use a change of variable to evaluate the integral $\int_0^x\mathrm{e}^{-\beta t}t^{2\nu}\,\mathrm{d}t=\frac{1}{\beta^{2\nu+1}}\gamma(2\nu+1,\beta x)$, which gives us inequality (\ref{ineqb2}).

\vspace{1mm}

\noindent{(ii)} Suppose now that $\nu\geq\frac{3}{2}$. A rearrangement of (\ref{55555}) gives that
\begin{align}&\int_0^x \mathrm{e}^{-\beta t}t^{\nu}\mathbf{L}_{\nu+1}(t)\,\mathrm{d}t-\beta \int_0^x \mathrm{e}^{-\beta t}t^{\nu}\mathbf{L}_{\nu}(t)\,\mathrm{d}t\nonumber\\
&\quad=\mathrm{e}^{-\beta x}x^\nu \mathbf{L}_\nu(x)-2\nu\int_0^x \mathrm{e}^{-\beta t}t^{\nu-1}\mathbf{L}_{\nu}(t)\,\mathrm{d}t-\int_0^x\mathrm{e}^{-\beta t}\frac{t^{2\nu}}{\sqrt{\pi}2^\nu\Gamma(\nu+\frac{3}{2})}\,\mathrm{d}t\nonumber\\
\label{1stint}&\quad=\mathrm{e}^{-\beta x}x^\nu \mathbf{L}_\nu(x)-2\nu\int_0^x \mathrm{e}^{-\beta t}t^{\nu-1}\mathbf{L}_{\nu}(t)\,\mathrm{d}t-\frac{\gamma(2\nu+1,\beta x)}{\sqrt{\pi}2^\nu\beta^{2\nu+1}\Gamma(\nu+\frac{3}{2})}.
\end{align}
We use inequality (\ref{Imon}) to bound the first integral on the left-hand side in (\ref{1stint}), and then divide through by $(1-\beta)$ and apply inequality (\ref{Imon}) again to obtain
\begin{align}\int_0^x \mathrm{e}^{-\beta t}t^{\nu}\mathbf{L}_{\nu}(t)\,\mathrm{d}t&\quad>\frac{1}{1-\beta}\bigg\{\mathrm{e}^{-\beta x}x^\nu \mathbf{L}_\nu(x)-2\nu\int_0^x \mathrm{e}^{-\beta t}t^{\nu-1}\mathbf{L}_{\nu}(t)\,\mathrm{d}t\nonumber\\
&\quad-\frac{\gamma(2\nu+1,\beta x)}{\sqrt{\pi}2^\nu\beta^{2\nu+1}\Gamma(\nu+\frac{3}{2})}\bigg\}\nonumber\\
&\quad>\frac{1}{1-\beta}\bigg\{\mathrm{e}^{-\beta x}x^\nu \mathbf{L}_\nu(x)-2\nu\int_0^x \mathrm{e}^{-\beta t}t^{\nu-1}\mathbf{L}_{\nu-1}(t)\,\mathrm{d}t\nonumber\\
\label{1stint0}&\quad-\frac{\gamma(2\nu+1,\beta x)}{\sqrt{\pi}2^\nu\beta^{2\nu+1}\Gamma(\nu+\frac{3}{2})}\bigg\}.
\end{align}
Lastly, we bound the integral $\int_0^x\mathrm{e}^{-\beta t}t^{\nu-1} \mathbf{L}_{\nu-1}(t)\,\mathrm{d}t$ using inequality (\ref{gau1}) (which can be done because $\nu\geq\frac{3}{2}$), which gives us inequality (\ref{ineqb3}).

\vspace{1mm}

\noindent{(iii)} Let $\nu>-1$, which ensures that all integrals in this proof of inequality (\ref{ineqb4}) exist.  We start with the same integration by parts to part (i), but with $\nu$ replaced by $\nu+1$:
\begin{align}\label{reart}\int_0^x \mathrm{e}^{-\beta t}t^{\nu+1}\mathbf{L}_{\nu+1}(t)\,\mathrm{d}t&=-\frac{1}{\beta}\mathrm{e}^{-\beta x}x^{\nu+1}\mathbf{L}_{\nu+1}(x)+\frac{1}{\beta}\int_0^x\mathrm{e}^{-\beta t}t^{\nu+1}\mathbf{L}_\nu(t)\,\mathrm{d}t,
\end{align}
where it should be noted that we used that $\lim_{x\downarrow0}x^{\nu+1}\mathbf{L}_{\nu+1}(x)=0$ for $\nu>-1$ (see (\ref{Itend0})). 
We now note the simple inequality $\int_0^x\mathrm{e}^{-\beta t}t^{\nu+1}\mathbf{L}_\nu(t)\,\mathrm{d}t<x\int_0^x\mathrm{e}^{-\beta t}t^{\nu}\mathbf{L}_\nu(t)\,\mathrm{d}t$, $x>0$, which holds because $\mathbf{L}_\nu(t)>0$ for $t>0$, $\nu>-1$.  Applying this inequality to (\ref{reart}) and rearranging gives
\begin{equation}\label{jj27}\int_0^x\mathrm{e}^{-\beta t}t^{\nu}\mathbf{L}_\nu(t)\,\mathrm{d}t>\mathrm{e}^{-\beta x}x^{\nu}\mathbf{L}_{\nu+1}(x)+\frac{\beta}{x}\int_0^x \mathrm{e}^{-\beta t}t^{\nu+1}\mathbf{L}_{\nu+1}(t)\,\mathrm{d}t.
\end{equation}
We can use (\ref{jj27}) to obtain another inequality
\begin{align*}&\int_0^x\mathrm{e}^{-\beta t}t^{\nu}\mathbf{L}_\nu(t)\,\mathrm{d}t\\
&\quad>\mathrm{e}^{-\beta x}x^{\nu}\mathbf{L}_{\nu+1}(x)+\frac{\beta}{x}\bigg(\mathrm{e}^{-\beta x}x^{\nu+1}\mathbf{L}_{\nu+2}(x)+\frac{\beta}{x}\int_0^x \mathrm{e}^{-\beta t}t^{\nu+2}\mathbf{L}_{\nu+2}(t)\,\mathrm{d}t\bigg)\\
&\quad=\mathrm{e}^{-\beta x}x^{\nu}\mathbf{L}_{\nu+1}(x)+\beta \mathrm{e}^{-\beta x}x^{\nu}\mathbf{L}_{\nu+2}(x)+\frac{\beta^2}{x^2}\int_0^x \mathrm{e}^{-\beta t}t^{\nu+2}\mathbf{L}_{\nu+2}(t)\,\mathrm{d}t,
\end{align*}
 and iterating gives inequality (\ref{ineqb4}).  In performing this iteration, it should be noted that the series $\sum_{k=0}^\infty \beta^k \mathbf{L}_{\nu+k+1}(x)$ is convergent.  This can be seen by applying inequality (\ref{Imon}) (since $\nu>-1$) to obtain that, for all $x>0$, $\sum_{k=0}^\infty \beta^k \mathbf{L}_{\nu+k+1}(x)<\mathbf{L}_{\nu+1}(x)\sum_{k=0}^\infty \beta^k=\frac{\mathbf{L}_{\nu+1}(x)}{1-\beta}$, with the assumption that $0<\beta<1$ ensuring that the geometric series is convergent.



\vspace{1mm}

\noindent{(iv)} Lastly, we prove that inequalities (\ref{ineqb2})--(\ref{ineqb4}) are tight in the limit $x\rightarrow\infty$. To this end, we note the following limiting forms, which hold for all $\nu>-1$ and $0<\beta<1$:
\begin{align}\label{eqeq1} \int_0^x \mathrm{e}^{-\beta t}t^\nu  \mathbf{L}_{\nu}(t)\,\mathrm{d}t&\sim \frac{1}{\sqrt{2\pi}(1-\beta)}x^{\nu-1/2}\mathrm{e}^{(1-\beta)x}, \quad x\rightarrow\infty,\\
\label{eqeq2}\mathrm{e}^{-\beta x}x^\nu \mathbf{L}_{\nu+n}(x)&\sim  \frac{1}{\sqrt{2\pi}}x^{\nu-1/2}\mathrm{e}^{(1-\beta)x}, \quad x\rightarrow\infty,\:n\in\mathbb{R},
\end{align}
where (\ref{eqeq2}) is immediate from (\ref{Itendinfinity}), and (\ref{eqeq1}) follows from using (\ref{Itendinfinity}) and a standard asymptotic analysis.
The tightness of inequalities (\ref{ineqb2}) and (\ref{ineqb3}) in the limit $x\rightarrow\infty$ follows immediately from  (\ref{eqeq1}) and (\ref{eqeq2}).  To show that inequality (\ref{ineqb4}) is tight as $x\rightarrow\infty$ we just need to additionally use that $\sum_{k=0}^\infty\beta^k=\frac{1}{1-\beta}$, since $0<\beta<1$. \hfill $\square$

\vspace{3mm}

\noindent{\emph{Proof of Theorem \ref{tiger2}.}} (i) Rearranging inequality (\ref{term00})  gives that, for $x>0$, $-\frac{1}{2}<\nu<\frac{1}{2}$, $0<\beta<1$,
\begin{align*}\int_0^x\mathrm{e}^{-\beta t}t^\nu \mathbf{L}_\nu(t)\,\mathrm{d}t<\frac{14}{(2\nu+1)(1-\beta)}\frac{\mathrm{e}^{-\beta x}x^{\nu-1}}{K_{\nu+3}(x)}.
\end{align*}
From the bound $\frac{1}{K_{\nu+3}(x)}<2x\mathbf{L}_{\nu+1}(x)$, which is a rearrangement of the lower bound in (\ref{klineq1}), we obtain that, for $x>0$, $-\frac{1}{2}<\nu<\frac{1}{2}$, $0<\beta<1$,
\begin{equation*}\int_0^x\mathrm{e}^{-\beta t}t^\nu \mathbf{L}_\nu(t)\,\mathrm{d}t<\frac{28}{(2\nu+1)(1-\beta)}\mathrm{e}^{-\beta x}x^\nu \mathbf{L}_{\nu+1}(x).
\end{equation*}
Using that $2\nu+1>0$ for $-\frac{1}{2}<\nu<\frac{1}{2}$ gives us
inequality (\ref{ineqb10}) for the case $-\frac{1}{2}<\nu<\frac{1}{2}$.  Inequality (\ref{ineqb10}) can in fact be seen to hold for all $\nu>-\frac{1}{2}$, by noting that the upper bound in inequality (\ref{ineqb10}) is strictly greater than the the upper bound in inequality (\ref{gau1}) (due to \cite{gaunt ineq4}), which is valid for $\nu\geq\frac{1}{2}$.

\vspace{1mm}

\noindent{(ii) We argue as in part (i), but we apply inequality (\ref{term10}), rather than inequality (\ref{term00}), and then use the bound $\frac{1}{K_{\nu+2}(x)}<2x\mathbf{L}_{\nu}(x)$.

\vspace{1mm}

\noindent{(iii) The proof proceeds exactly as that of inequality (\ref{ineqb3}), with the sole modification being that we use (\ref{ineqb10}) to bound the integral on the right-hand side of (\ref{1stint0}), instead of inequality (\ref{gau1}).   The tightness of inequality (\ref{ineqb12}) in the limit $x\rightarrow\infty$ is established by the same argument as that used in part (iv) of the proof of Theorem \ref{tiger1}.
\hfill $\square$

\vspace{3mm}

\noindent{\emph{Proof of Proposition \ref{propone}.}} (i) To get inequality (\ref{ineqb21}), in part (i) of the proof of Theorem \ref{tiger1} use inequality (\ref{Imon}) to bound the third integral in (\ref{55555}), instead of the first integral.

\vspace{1mm}

\noindent{(ii)} To get inequality (\ref{ineqb22}), in part (ii) of the proof of Theorem \ref{tiger1} use (\ref{Imon}) to bound the second integral in (\ref{1stint}), instead of the first integral.

\vspace{1mm}

\noindent{(iii)} By studying the proof of inequality (\ref{ineqb12}), it can be seen that the above alteration that gave us inequality (\ref{ineqb22}) instead of inequality (\ref{ineqb3}) can also be used to give us inequality (\ref{ineqb23}). \hfill $\square$

\appendix

\section{Basic properties of modified Struve and modified Bessel functions}\label{appa}

In this appendix, we present some basic properties of the modified Struve function of the first kind $\mathbf{L}_\nu(x)$ and the modified Bessel functions $I_\nu(x)$ and $K_\nu(x)$ that are used in this paper. All formulas are given in \cite{olver}, except for the inequality which was obtained by \cite{bp14}.  

The modified Struve function $\mathbf{L}_\nu(x)$ is a regular function of $x\in\mathbb{R}$, and is positive for all $\nu\geq-\frac{3}{2}$ and $x>0$. The modified Bessel functions $I_{\nu}(x)$ and $K_{\nu}(x)$ are also both regular functions of $x\in\mathbb{R}$.  For $x>0$, the functions $I_{\nu}(x)$ and $K_{\nu}(x)$ are positive for $\nu\geq-1$ and all $\nu\in\mathbb{R}$, respectively.
The modified Struve function $\mathbf{L}_\nu(x)$ satisfies the following recurrence relation and differentiation formula
\begin{align}\label{Iidentity}\mathbf{L}_{\nu -1} (x)- \mathbf{L}_{\nu +1} (x) &= \frac{2\nu}{x} \mathbf{L}_{\nu} (x)+\frac{(\frac{1}{2}x)^\nu}{\sqrt{\pi}\Gamma(\nu+\frac{3}{2})}, \\
\label{diffone}\frac{\mathrm{d}}{\mathrm{d}x} \big(x^{\nu} \mathbf{L}_{\nu} (x) \big) &= x^{\nu} \mathbf{L}_{\nu -1} (x).
\end{align}
We have the following special cases
\begin{equation}\label{speccase}\mathbf{L}_{-\frac{1}{2}}(x)=\sqrt{\frac{2}{\pi x}}\sinh(x),\quad I_{-\frac{1}{2}}(x)=\sqrt{\frac{2}{\pi x}}\cosh(x).
\end{equation}
We also have the following asymptotic properties:
\begin{align}\label{Itend0}\mathbf{L}_{\nu}(x)&\sim \frac{x^{\nu+1}}{\sqrt{\pi}2^\nu\Gamma(\nu+\frac{3}{2})}\bigg(1+\frac{x^2}{3(2\nu+3)}\bigg), \quad x \downarrow 0, \: \nu>-\tfrac{3}{2}, \\
\label{Itendinfinity}\mathbf{L}_{\nu}(x)&\sim \frac{\mathrm{e}^{x}}{\sqrt{2\pi x}}\bigg(1-\frac{4\nu^2-1}{8x}\bigg), \quad x \rightarrow\infty, \: \nu\in\mathbb{R}, \\
\label{Ktend0}K_{\nu} (x) &\sim \frac{2^{\nu-1}\Gamma(\nu)}{x^\nu},\quad \:x\downarrow0,\: \nu>0, \\
\label{Ktendinfinity} K_{\nu} (x) &\sim \sqrt{\frac{\pi}{2x}} \mathrm{e}^{-x}\bigg(1+\frac{4\nu^2-1}{8x}\bigg), \quad x \rightarrow \infty, \: \nu\in\mathbb{R}.
\end{align}
It was shown by \cite{bp14} that, for $x>0$,
\begin{equation}\label{Imon}\mathbf{L}_{\nu} (x) < \mathbf{L}_{\nu - 1} (x), \quad \nu \geq \tfrac{1}{2}.  
\end{equation} 
Other inequalities for the modified Struve function $\mathbf{L}_\nu(x)$ are given in \cite{bp14,bps17,gaunt ineq5,jn98}, some of which improve on inequality (\ref{Imon}).

\subsection*{Acknowledgements}
I would like to thank the referees for their helpful comments and suggestions that helped me improve the presentation of my paper.

\bibliographystyle{amsplain}




\end{document}